\documentclass{conm-p-l}

\usepackage[shortlabels]{enumitem}
\usepackage{aliascnt}

\usepackage{amsmath,amssymb, amsfonts}
\usepackage{amsthm}
\usepackage{hyperref}

\hypersetup{
  colorlinks   = true, %Colours links instead of ugly boxes
  urlcolor     = blue, %Colour for external hyperlinks
  linkcolor    = blue, %Colour of internal links
  citecolor   = red %Colour of citations
}

\usepackage[msc-links]{amsrefs}

\usepackage[capitalise]{cleveref}

\crefname{thm}{Thm.}{}
\crefname{prop}{Prop.}{}
\crefname{lem}{Lem.}{}
\crefname{cor}{Cor.}{}

\newtheorem*{thm*}{Theorem}
\newtheorem{thm}{Theorem}
\newtheorem{prop}{Proposition}
\newtheorem{lem}{Lemma}
\newtheorem{cor}{Corollary}
\newtheorem{rem}{Remark}
\newtheorem{exa}{Example}
\newtheorem{prob}{Problem}

\theoremstyle{definition}
\newtheorem{defi}{Definition}

\theoremstyle{remark}

\DeclareMathOperator\re{Re }

\def\cR{\mathcal R}

\def\Z{\mathbb Z}
\def\R{\mathbb R}
\def\Q{\mathbb Q}
\def\C{\mathbb C}
\def\P{\mathbb P}
\def\wP{\mathbb{WP}}

\def\C{{\mathbb C}}

\def\H{\mathcal H}

\def\M{\mathcal M}
\def\X{\mathcal X}

\def\O{\mathcal O}
\def\p{\mathfrak p}

\def\A{\mathbb A}

\def\X{\mathcal X}
\def\Y{\mathcal Y}

\def\a{{\alpha }}
\def\b{{\beta }}

\def\D{\Delta}

\def\Aut{\mbox{Aut }}

\DeclareMathOperator\twist{\mbox{Twist }}

\DeclareMathOperator\Gal{Gal}

\DeclareMathOperator\Orb{Orb}
\DeclareMathOperator\homo{Hom}

\def\l{\lambda}
\def\a{\alpha}

\def\b{\beta}

\newcommand\x{\mathbf{x}}

\newcommand\w{\mathbf{w}}

\def\I{\mathcal I}

\DeclareMathOperator\wgcd{\mathit{wgcd }}   %   weighted gcd
\DeclareMathOperator\wh{\mathfrak{h}}   %  weighted height

\DeclareMathOperator\awh{\mathfrak{\tilde h} }   % absolute    weighted height

\def\awgcd{\overline{\wgcd }}

\def\q{\mathfrak q}
\def\J{\mathfrak J}

\def\D{\Delta}

\def\Aut{\mbox{Aut }}

\def\<{\left<}
\def\>{\right>}

\def\val{\mathbf{val}}

\def\RR{\mathcal R}

\def\q{\mathfrak q}

\DeclareMathOperator\SL{SL}
\DeclareMathOperator\GL{GL}

\DeclareMathOperator\h{H}

\newcommand\Df{\mathfrak D_f}

\newcommand\z{\xi}  %zero map
\newcommand\F{\mathcal F}            % Fundamental domain
\def\G{\Gamma}

\title{Reduction of superelliptic Riemann surfaces}
\author{Tanush Shaska}
\address{Research Institute of Science and Technology (RISAT),
Vlor\"{e}, ALBANIA}
\email{shaska@risat.org}

\begin{document}

\begin{abstract}
For a superelliptic curve $\mathcal X$, defined over $\Q$, let $\mathfrak p$ denote the corresponding moduli point in the weighted moduli space. 
We describe a method how  to determine  a minimal integral model of  $\X$ such that: 
i) the corresponding moduli point $\p$  has minimal weighted height,  ii)   the equation of the curve has minimal coefficients.   
Part i) is accomplished by reduction of the moduli point which is equivalent with obtaining a representation of the moduli point $\p$ with minimal weighted height,   as defined in \cite{b-g-sh}, and part ii) by the classical reduction of the binary forms.   
\end{abstract}

\maketitle              
%\setcounter{tocdepth}{2}
%\tableofcontents
%****************************************
\section{Introduction}
Let $k$ be an algebraic number field   and  $\O_k$ its ring of integers.  The isomorphism class of a  smooth, irreducible, planar, algebraic curve $\X$, defined over $\O_k$, is determined by its set of invariants which are  homogenous polynomials in   coefficients of $\X$.  The best understood case is when $\X$ is a hyperelliptic.  In    \cite{m-sh, frey-sh} the authors make the case that superelliptic curves are a natural generalization of hyperelliptic curves and asked if many questions and arithmetic results of   hyperelliptic curves can be extended to superelliptic curves.  
%While questions asked in \cite{m-sh} cover a wide variety of problems, 
Here we focus on minimal models for the curve and minimal representation for the corresponding moduli point.  

By a \emph{superelliptic curve} we mean  a smooth, irreducible, plannar,  algebraic curve $\X$, defined over $k$,  with projective   equation $z^m y^{d-m} = f(x, y)$, where $f(x, y)$ is a degree $d$  binary form of nonzero discriminant $\D_f\neq 0$. We assume that such curves have a normal equation over $\O_k$, in other words $f(x, y) \in \O_k[x, y]$; see \cite{m-sh} for details. The isomorphism class of  such curves is determined by the set of invariants of the degree $d$ binary forms. These invariants are generators of the invariant ring of binary forms of fixed degree; see \cite{rota}  among other places.
%Such invariant rings are theoretically well understood and their generators have been known explicitly for degree $d \leq 8$ due to work of classical invariant theorists such as Clebsch, Bolza, Gordan, van Gall, et al.   In the last two decades computational invariant theory has been revived and have been attempts to compute such generators for higher degree binary forms. 

Given a binary form $f(x, y)$, there are two main \emph{reduction} problems. Determine   $g(x, y)$,  $\GL_2(k)$-equivalent to $f(x, y)$, such that $g(x, y)$ has minimal:
\begin{itemize}[\upshape(a), nolistsep]
\item[a)]   invariants
\item[b)]  coefficients
\end{itemize}
Explaining  what \emph{minimal} means  will  be the main focus of this paper.    Completing of each task we will call \textbf{reduction type a)} or \textbf{moduli reduction}  and  \textbf{reduction type b)} or \textbf{classical reduction}.    
For the purposes of this paper we will focus on the case when $k=\Q$, even though many of the results follow for any number field $k$. 
%
\iffalse
  Throughout this paper the following notation will be used:

%\noindent \textbf{Notation:} Throughout this paper the following notation will be used:

$k$  an algebraic number field

$\O_k$  the ring of integers of $k$

$\X$ a superelliptic curve with projective equation $z^m y^{d-m} = f(x, y)$

$f(x, y)$ a degree $d\geq 5$ binary form

$\D_f$ the discriminant of $f(x, y)$

$I_0, \ldots , I_n$  generators of the ring of invariants of degree $d$ binary forms 

\fi
%Reduction a) described in this paper, to the  best of our knowledge is new and has not appeared in the literature before, apart from elliptic curves in the work of Birch and Swinnerton-Dyer as explained below  (cf. \cref{section-4-3}).    

%Reduction b) is explained in details in \cref{section-4-2} and is based on work of G. Julia  \cite{julia}. 

%*****************************************************
\subsection{Reduction a): Minimality of the moduli point}
The   isomorphism class of $\X$ correspond to the equivalence class of  a  binary form $f(x, y)$. Thus, the isomorphism classes of superelliptic curves over $\bar k$) are determined by the set of generators of the ring of invariants $\cR_d$ of degree $d\geq 2$ binary forms,  in other words by $\SL_2 (\O_k)$-invariants.  By Hilbert's basis theorem $\cR_d$ is finitely generated. 
Let $I_0$, $\dots $, $I_n$ be the generators of   $\cR_d$   such that the homogenous degree of $I_i$ is $q_i$ for each $i=0, \dots , n$.  We denote by $\I:= (I_0, \dots , I_n)$ the tuple of invariants and by $\I(f)= \left(I_0 (f), \dots , I_n (f)\right)$ such invariants evaluated at the   binary form  $f$. The corresponding set of invariants  $\I(f)$ determines a point $\p= [\I(f)]$ in the weighted projective space $\wP_\w^n (k)$ (cf. \cref{w-space}) with weights $\w=(q_0, \dots , q_n)$; see  \cite{b-g-sh, m-sh} for details.  
 Denote the weighted  greatest common divisor of the weighted tuple  $\I(f)$ by $\wgcd (\I(f))$   and by $\awgcd (\I(f))$ its absolute weighted  greatest common divisor (cf. \cref{w-gcd}. 
The curve $\X$ is said to have a \textbf{minimal modular model} over $\O_k$ when it has minimal height   (cf. \cref{def:height}) in the weighted moduli space $\wP_\w^n (k)$. 
In other words, when the weighted valuation of the tuple $\I(f)$ for each prime $p \in \O_k$, 
$ \val_p  ( \I (f)) := \max   \left\{  \nu_p (I_i (f)) \; \text{ for all } \; i=0, \dots, n     \right\}  $
is minimal, where $\nu_p (I_i (f))$ is the valuation at $p$ of $I_i (f)$. The weighted moduli point $\p=[\I(f)]$ is called   \textbf{normalized} when $\wgcd (\I(f))=1$.  
%A given curve $\X$ which has normalized moduli point (cf. \cref{normalized})   is already in the minimal model. 
There are two main tasks:
\begin{enumerate}[\upshape(i), nolistsep]
\item   determine an equation of a curve $\X^\prime$, defined over $\O_k$, and $k$-isomorphic to $\X$,  with \emph{minimal modular model}. 
\item  determine an equation of a curve $\X^{\prime\prime}$, defined over $\O_k$ and $\bar k$-isomorphic to $\X$ with minimal modular model, in other words a twist of $\X$ with minimal modular model.
\end{enumerate}
Our main result is  \cref{thm-1}, which says     that  a minimal modular model for superelliptic curves exist.     Moreover, an  equation  $\X :   z^m y^{d-m}= f(x, y)$   is a minimal modular model over $\O_k$,   if   for every prime $p \in \O_k$ which divides  $ p \, \mid \,  \wgcd \left( \I(f) \right)$,  the valuation $\val_p $  of $\I(f)$ at $p$ satisfies  $\val_p (\I(f))  <   \frac d 2  \,  q_i$,    for all $i=0, \dots , n$.  Additionally,   for  $\lambda=\wgcd (\I (f)) $ with respect  the weights   $\left(\left\lfloor  \frac {dq_0} {2} \right\rfloor,\ldots, \left\lfloor  \frac {dq_n} {2} \right\rfloor \right) $
the transformation $ (x, y, z) \to \left(      \frac x \lambda, y, \lambda^{\frac d m} z \right) $    gives a minimal model of  $\X$ over $\O_k$. If $m|d$ then this isomorphism is defined over $k$. 

Such minimal modular models can be found not only up to $k$-isomorphism, but also over its algebraic closure $\bar k$. We call them \emph{minimal modular twists}.  In \cref{thm-2} we prove that  minimal modular twists of  superelliptic curves exist.   An   equation $\X : \;  z^m y^{d-m}= f(x, y)$ is a minimal modular twist  over $\O_k$,   if   for every prime $p \in \O_k$ such that   $ p \, \mid \,  \awgcd \left( \I(f) \right)$,  the valuation $\val_p $  of $\I(f)$ at $p$ satisfies  $\val_p (\I(f))  <   \frac d 2  \,  q_i$,   for all $i=0, \dots , n$. 

 \cref{thm-1}  makes it possible to create a database of superelliptic curves defined over $\O_k$, by storing only the curves with minimal weighted moduli point     and \cref{thm-2}  does this also for all their twists as well.    Our  motivation of exploring such reduction came from   genus 2 curves; see \cite{rat-pts}. 
%The idea of weighted heights and weighted greatest common divisors were developed in \cite{b-g-sh} and the moduli reduction was first developed in \cite{m-sh}.

%***************
\subsection{Reduction b):   Minimality of coefficients}
Let $\p \in \wP_\w^n (\Q)$ be a    normalized  moduli point via \emph{reduction a)} and $\X$   the corresponding superelliptic curve.
Assume that $\X$    has   equation $z^m y^{d-m} = f(x, y)$ over some minimal field of definition $k$.  Determining such equation is part of the math folklore for elliptic curves.  
%For genus 2 such equation is determined in  \cite{mestre} for  curves with automorphism group of order 2   and in  \cite[Theorem~3]{sh-com-alg} for all groups of order $> 2$.  
%Recently in  \cite{genus-2-univ}  a universal equation for all genus 2 curves is given. 
%Similar work has been done for higher genus for hyperelliptic curves; see \cite{g-sh} among others. Equations of cyclic curves with extra automorphisms have been determined in \cite{san-sh}.
%
One of the main concerns when obtaining such equations is to have coefficients as small as possible.   This comes down to  what historically is referred to as reduction of binary forms and goes back to Hermite for quadratic forms and Julia \cite{julia} for cubic forms. In \cite{stoll-cremona} such reduction was considered again for cubics and quartics and in \cite{beshaj, beshaj-thesis, beshaj-2} for $d\geq 3$.

Let $f(x, y)$ be a degree $d\geq 2$ binary form with real coefficients.  Then its \textbf{Julia quadratic} $\J_f$  is defined as in \cref{Julia-quadratic}. It is a positive definite quadratic and therefore has one root  in the upper-half complex plane $\H_2$, say $\alpha_f$.  Since $\J (f)$ is an $\SL_2 (\Z)$-covariant,  then bringing $\a_f$ to the fundamental domain  $\F$  by a matrix $M\in \SL_2(\Z)$, induces an action $f \to f^M$ on binary forms. The form $f^M$ is called \textbf{reduction of $f$}.   In \cite{beshaj} is given an approach of how to determine a minimal model for any degree binary form $f(x, y)$ and a minimal twist of $f(x, y)$ (cf. \cref{min-model}). 

%**************** 
\subsection{Birch and Swinnerton-Dyer computations}
After the first draft of this paper was written we discovered that a reduction combining both methods described above had been used in the seminal paper of Birch and Swinnerton-Dyer in \cite{birch-s-d-1} and \cite{birch-s-d-2} for computations with elliptic curves.  While the case of the elliptic curves is simpler, it is also the only case that is fully understood, since only for cubics and quartics we have precise results of the Julia reduction \cite{julia} as described in the work of Cremona and Stoll in \cite{stoll-cremona} and   by Beshaj in   \cite{beshaj, beshaj-thesis, beshaj-2} as remarked above.  

The initial computations in \cite{birch-s-d-1} start with reduction a).  
Let $f(x,y)=ax^4+bx^3 + cx^2+dx+e$, where $a, b, c, d, e$ are rational integers,  and $I_4, I_6$ its invariants.  If $p \in \Z$ is a prime such that $p\neq 2, 3$,  and $p | \wgcd_{4, 6} (I_4, I_6)$, there there a quartic integral binary form $g(x, y)$, $\GL_2 (\Z)$ equivalent to $f(x, y)$ with invariants $(p^{-4} I_4, p^{-6} I_6)$.  Birch and Swinnerton-Dyer go to a case by case analysis to prove this result; see \cite[Lemma~3]{birch-s-d-1}. However, this is an immediate consequence of the reduction a) described in \cref{red-a}.  Lemma~4 and Lemma~5 in \cite{birch-s-d-1} give the reduction for $p=2$ and 3 respectively. 

Minimal integral models   give the "nicest" equation of the paper over a global field, since the corresponding invariants are non-zero for as many primes $p \in \O_k$ as possible.  It is a topic of interest to explore the stability of such curves.  Adjusting the method of reduction a) would give a method of obtaining a stable or semistable form for any binary form; see \cite{curri}. 

Both reductions can be performed for all superelliptic curves, providing that we explicitly know the generators of the ring of invariants $\cR_d$.  For non-superelliptic curves this is a much more difficult problem which requires the full arsenal of GIT.  An explicit  description of the moduli point is not known in general.    Moreover, there is no known algorithm to determine the equation of the curve starting with the  moduli point even in the case a planar curves.  
%Superelliptic curves are a natural extension of the theory of hyperelliptic curves as explored in \cite{m-sh} and \cite{book}.  Extending this approach   beyond superelliptic curves seems very difficult at the moment. 

\medskip

\noindent \textbf{Acknowledgments:}   I would like to thank Mike Fried for helpful comments and discussions.

%*****************************************
\section{Preliminaries}
\subsection{Superelliptic curves}  
A  a genus $g \geq 2$ smooth, irreducible,  algebraic curve    $\X$ defined over an algebraically closed field $k$ is called  a \textbf{superelliptic curve of level $n$}  if there exist an element $\tau \in \Aut (\X)$ of order $n$ such that $\tau$ is central and  the quotient $\X / \< \tau \> $ has genus zero; see \cite{m-sh} for details.    Superelliptic curves  have affine equation 
\begin{equation}\label{super}
\X  : \; y^n = f(x) = \prod_{i=1}^d (x-\a_i), \quad \text{ for} \quad \D_f\neq 0.
\end{equation}

Denote by $\sigma : \X \to \X$ the superelliptic automorphism, i.e. $ \sigma (x, y) \to (x, \xi_n y)$,   where $\xi_n$ is a primitive $n$-th root of unity.  Notice that $\sigma$ fixes 0 and  the point at infinity in $\P_y^1$.    The natural projection  $ \pi : \X \to \P^1_x=\X/\<\sigma\>$    is called the \textbf{superelliptic projection}.  It  has  $\deg \pi =n$ and  $ \pi (x, y) =  x$.    This cover is branched at exactly at the roots $\a_1, \dots , \a_d$ of $f(x)$.   Then the affine equation is   $\X : \; z^m = \prod_{i=1}^d (x-\a_i)$.   Denote the projective equation of $\X$ by 
\begin{equation}\label{w-eq-super}
z^m y^{d-m}= f(x, y) = a_d x^d + a_{d-1}x^{d-1}y + \cdots + a_1 x y^{d-1} + a_0y^d
\end{equation}
defined over a field $k$.   Hence, $f(x, y)$ is a binary form of degree $\deg f =d$.   

Let  $k[x, y]$  be the  polynomial ring in  two variables and   $V_d (k)$ denote  the $(d+1)$-dimensional  subspace  of  $k[x, y]$  consisting of homogeneous polynomials  $f(x, y)$ of  degree $d$. Elements  in $V_d$  are called  \textit{binary  forms} of degree $d$.   The general linear group  $\GL_2(k)$ acts as a group of automorphisms on $k[x, y] $.

Consider $a_0$, $a_1$, $\ldots$ , $a_d$ as parameters  (coordinate  functions on $V_d$). Then the coordinate  ring of $V_d$ can be identified with $ k[a_0  , ... , a_d] $. For $I \in k[a_0, \ldots , a_d]$ and $M \in \GL_2(k)$, define $I^M \in k[a_0, \dots , a_d]$ as  ${I^M}(f):= I( f^M)$,  for all $f \in V_d$. Then  $I^{MN} = (I^{M})^{N}$ and we have  an action of $\GL_2(k)$ on $k[a_0, \ldots , a_d]$ (cf. \cref{GL-action}).
A homogeneous polynomial $I\in k[a_0, \dots , a_d, x, y]$ is called a \textbf{covariant}  of index $s$ if  $I^M(f)=\l^s I(f)$,  where $\l =\det  (M)^d$.  The homogeneous degree in $a_0, \dots , a_d$ is called the \textbf{ degree} of $I$,  and the homogeneous degree in $x, y$ is called the \textbf{  order} of $I$.  A covariant of order zero is called \textbf{invariant}.  
From Hilbert's basis theorem  the ring of invariants $\RR_d$  of degree $d$ binary forms is finitely generated. 
Let $I_0, \dots  , I_n$ be the generators of $\RR_d$ with degrees $q_0, \dots , q_n$ respectively.     Usually we assume an ordering  $q_0 < q_1 < \ldots < q_n$ of degrees.  Denote the ordered  tuple of invariants by $\I := (I_0, \dots , I_n)$.  
Over the algebraic closure  we have 
\begin{equation} f(x, y)= (y_1 x - x_1 y) \cdots (y_dx - x_d y)= \prod_{i=1}^d \det \begin{pmatrix} x & x_i \\ y & y_i \end{pmatrix}, \end{equation}
where $(x_i, y_i)$ are the homogenous coordinates of the roots in $\P^1$. 
For $ M  \in \GL_2(k) $, denote by $\delta = \det M$ and by  $f^M$ the action of $M$ on $f$.   Each root $(x_i, y_i)$ goes to  $M \begin{bmatrix} x_i \\ y_i\end{bmatrix}$ and 
$\begin{pmatrix} x & x_i \\ y & y_i \end{pmatrix} \to M \begin{pmatrix} x & x_i \\ y & y_i \end{pmatrix}$.  Hence, we have 
\begin{equation}\label{GL-action}
f^M (x, y) = \prod_{i=1}^d \det  \left( M \cdot    \begin{pmatrix} x & x_i \\ y & y_i \end{pmatrix} \right) = \left( \det M \right)^d \, f(x, y). 
\end{equation}    
Hence, all the coefficients $a_i$, $i=0, \ldots , d$ are multiplied by $\delta^d$.  Since an invariant of degree $s$ is a homogenous polynomial of degree $s$ in terms of $a_i$, 
\begin{equation} 
I_s = \sum  a_0^{\a_0} \ldots a_d^{\a_d}
\end{equation}
where $\a_i=0, \ldots , d$ and $\a_0+ \cdots + \a_d=s$,   then  $I (f^M) = \delta^{ds} \, I_s (f)$.    However, it is more difficult if we want to determine $f^M (x, y)$ when $f(x, y)$ is given as in \cref{w-eq-super}. Such expressions of invariants in terms of coefficients are determined via \emph{transvections} and \emph{umbral calculus}; see \cite{rota} among many other sources. 

Consider $f(x, y) $ as in  \cref{w-eq-super}, $M  =  \begin{bmatrix} \alpha & \beta \\ \gamma & \delta    \end{bmatrix}    \in \GL_2 (k)$, and $f^M (x, y):= f \left(  \alpha x + \beta y, \gamma x + \delta y     \right)$.  We want to determine the invariants of $f^M (x, y)$ in  terms of the invariants of $f(x, y)$. 
The following result is fundamental to our approach. 
\begin{prop}\label{prop-2}    
Let  $f  \in V_d (k)$, $M\in \GL_2 (k)$, and  $ I_0,  \ldots  I_n$ be the generators of $\RR_d$ with degrees $q_0, \ldots , q_n$ respectively. 
Then    for each $i=0, \ldots , n$ 
\begin{equation}\label{prop-2}
I_i  \left( f^M (x, y) \right)  =   \l^{q_i}\,  I_i (f)
%\left( I_0 (f), \dots   I_i (f), \dots , I_n (f) \right) = \left( \l^{q_0} \, I_0 (g), \dots ,    \l^{q_i}\,  I_i (g), \dots , \l^{q_n} \, I_n (g)   \right), 
\end{equation}
where   $\l = \left( \det M \right)^{\frac d 2}$. 
\end{prop}

\proof  Let $J$ be a covariant of degree $\deg J = s$, order $\mbox{ord } (J)$, and weight $\mbox{wt } (J)$.  Recall that the \textbf{weight} of a covariant is the integer $r$ such that 
\begin{equation}J \left( f \left(  \alpha x + \beta y, \gamma x + \delta y     \right) \right) = (\alpha \delta - \beta \gamma)^r J (f(x, y)).\end{equation}
We know that 
\begin{equation} \deg J + 2 \mbox{wt } (J) = (  \deg f + 2 \mbox{wt } f) \,    \mbox{ord} (J), \end{equation}
see \cite[Prop.~2.29]{olver}.   Let $I$ be an invariant of order $\mbox{ord } (I):=s$.  Then    $\deg I =0$,   $\mbox{wt } (f)=0$ and we have $2 \mbox{wt } (I)= d s$.  So the weight of any invariant is $\frac d 2 s$.  This completes the proof. 
\qed

%\begin{rem} 
%The above result was stated in   \cite[Prop.~14]{m-sh}, but the proof accidentally cut off.   It was probably known to classical invariant theorists, but we couldn't find it stated anywhere in the literature.  Determining weights of $\SL_2 (k)$ invariants has appeared in the literature mostly in terms of $\SL_2 (\R)$ heat kernels; see \cite{2020-2} among many others. 
%\end{rem}
%**********************************************************************
\subsection{$\mbox{Proj } \RR_d$ as a weighted projective space} \label{w-space}
Since all $I_0, \dots , I_i, \dots , I_n$ are homogenous polynomials then $\RR_d$ is a graded ring and  $\mbox{Proj } \RR_d$ is a weighted projective space. 
Let   $\w:=(q_0, \dots , q_n) \in \Z^{n+1}$ be the fixed ordered tuple of positive integers called \textbf{weights}.   Consider the action of $k^\star = k \setminus \{0\}$ on $\A^{n+1} (k)$ as follows
$ \lambda \star (x_0, \dots , x_n) = \left( \l^{q_0} x_0, \dots , \l^{q_n} x_n   \right) $, 
for $\l\in k^\ast$.  The quotient of this action is called a \textbf{weighted projective space} and denoted by   $\wP^n_{\w} (k)$. 
%The space $\P(1, \dots , 1)$ is the usual projective space. 
It is the projective variety $Proj \left( k [x_0,...,x_n] \right)$ associated to the graded ring $k [x_0, \dots ,x_n]$ where the variable $x_i$ has degree $q_i$ for $i=0, \dots , n$. 
We will denote a point $\p \in \wP_w^n (k)$ by $\p = [ x_0 : x_1 : \dots : x_n]$.    For proofs of the following two results see \cite{m-sh}. 

\begin{prop} Let $I_0,  I_1,   \dots,  I_n$ be the generators of the ring of invariants $\RR_d$ of degree $d$ binary forms. 
A $k$-isomorphism class of a binary form $f$ is determined by  the weighted moduli point 
\begin{equation} 
\I (f) := \left[ I_0 (f), I_1 (f), \dots , I_n(f) \right]  \in \wP_\w^n  (k). 
\end{equation}
Moreover, $f=g^M$ for some $M\in \GL_2 (K)$ if and only if  $ \I(f)  = \lambda\star \I(g)$, 
for $\lambda = \left( \det A \right)^{\frac d 2}$.
\end{prop}

Since the isomorphism class of a superelliptic curve    $\X : z^m y^{d-m}= f(x, y)$ is determined by the equivalence  class of binary form $f(x, y)$ we denote the set of invariants of $\X$ by $\I (\X) : = \I (f)$. 
\begin{cor}
Let $\X$ be  as in \cref{w-eq-super}. The $\bar k$-isomorphism class of $\X$ is determined by the weighted moduli point  $\p := \left[ \I (f) \right]  \in \wP_\w^n  (k)$. 
\end{cor}

%********************************************************************************
\subsection{Heights}\label{heights}
If we want to create a list of all isomorphism classes of superelliptic curves we have to create a database of points in $\wP_\w^n (k)$.  Obviously,  we would prefer to take for each point $\p \in \wP_\w^n (k)$ its \emph{smallest} representative and  also want to order such points.  We struggled with both such tasks in \cite{rat-pts}, until we were able to define a height function on $\wP_\w^n (k)$, which solves both problems.  This is due to a Northcott like theorem for this weighted height, which says that there are only finitely many points for a bounded height.  
Also the representation of the point $\p \in \wP_\w^n (k)$ with \emph{smallest coordinates} correspond precisely to the tuple $(x_0, \ldots , x_n)$ such that its weighted greatest common divisor is $=1$.  These ideas were developed in detail in \cite{b-g-sh}.

Fix the following notation: $k$  is a number field,   $\O_k$ is ring of integers of $k$,   $M_k$ is a complete set of absolute values of $k$,  $M_k^0$ is   the set of all non-archimedian places in $M_k$,  $M_k^\infty$ is   the set of Archimedean places, and   $\X/k$ is   a smooth projective superelliptic curve defined over $k$.     
For a place $\nu \in M_k$, the corresponding absolute value is denoted by $| \cdot |_\nu$, normalized with respect to $k$ such that the product formula holds and the Weil height for $x \in k$  is 
$H (x) =  \prod_\nu \max \{ 1, |x|_\nu \}$.

For  $P=[x_0, \dots , x_n] \in \P^n(k)$ the \textbf{multiplicative heigh}t of $P$ is defined as follows
\begin{equation}
H_k (P) := \prod_{v \in M_k} \max\left\{\frac{}{}|x_0|_v^{n_v} , \dots, |x_n|_v^{n_v}\right \},
\end{equation}
where $n_v$ is the \textbf{local degree at $v$} given by $n_v= [ k_v : \Q_v]$ for $k_v$ and $\Q_v$ are the completions with respect to $v$. 
The height $H_k (P)$ is well defined, in other words it does not depend on the choice of homogenous coordinates of $P$.   Moreover,  $H_k (P) \geq 1$.
 If $L/k$ is a finite extension, then    $H_L(P)=H_k (P)^{[L:k]}$.  
Hence,  we can define the height on $\P^n (\overline \Q)$, which  is called the \textbf{absolute (multiplicative) height} and is the function 
$H: \P^n (\bar \Q)  \to [1, \infty)$, such that $H(P)=H_k(P)^{1/[k:\Q]}$.  
The height is invariant under Galois conjugation. In other words,  if  $P \in \P^n(\overline \Q)$ and $\sigma \in G_{ \Q}$, then    $H(P^\sigma) = H(P)$.
The following are the two main results in the theory of heights on projective spaces; see \cite{bombieri} or \cite{silv-book}. 

\begin{thm}[Northcott]  \label{thm_finite}
Let $c_0$ and $d_0$ be constants. Then the set 
\begin{equation}\{P \in \P^n(\overline \Q): H(P) \leq c_0 \text{ and } [\Q(P):\Q] \leq d_0\}\end{equation}
has  finitely many points. In particular,   $\{P \in \P^n(k): H_k(P) \leq c_0 \}$  is a finite set.
\end{thm}

\begin{thm}[Kronecker's theorem]  
Let 
%$k$ be a number field,  and let 
$P=[x_0, \dots , x_n] \in \P^n(k)$. Fix any $i_0$ with $x_{i_0} \neq 0$. Then $H (P)= 1$ if and only if the ratio $x_j/x_{i_0}$ is a root of unity or zero for every $0 \leq j \leq n$. 
\end{thm}

%*******************
%\subsubsection{Heights on weighted projective spaces}

Amazingly, the above theory can be extended to weight projective spaces with necessary adjustments. We give a quick recapture here; see \cite{b-g-sh} for details.  
Let $\x = (x_0, \dots x_n ) \in \Z^{n+1}$ be a tuple of integers, not all equal to zero.  
 A \textit{weighted integer tuple} is a tuple $\x = (x_0, \dots, x_n ) \in \Z^{n+1}$ such that to each coordinate $x_i$ is assigned the weight $q_i$ and  $\w:=(q_0, \dots , q_n)$ is called the  \emph{set of weights}. 
 We multiply weighted tuples by scalars $\lambda \in \Q$ via 
\begin{equation} \lambda \star (x_0, \dots , x_n) = \left( \l^{q_0} x_0, \dots , \l^{q_n} x_n   \right) \end{equation}
For an ordered   tuple of integers  $\x=(x_0, \dots, x_n) \in \Z^{n+1}$, whose coordinates are not all zero, the \textbf{weighted greatest common divisor with respect to the set of weights} $\w$ is the largest integer $d$ such that 
$ d^{q_i} \, \mid \, x_i, \; \; \text{for all }    i=0, \dots, n$. 
%see \cite{mandili, b-g-sh} for details.    
For a weighted tuple  $\x =(x_0,\ldots,x_n) \in \O_k^{n+1}$ the weighted greatest common divisor is given by 
\begin{equation}\label{w-gcd}
\wgcd (\x) = \prod_{p \in \O_k}    p^{ \min \left\{ \left\lfloor  \frac {\nu_p (x_0)}  {q_0}   \right\rfloor,  \dots ,  \left\lfloor  \frac {\nu_p (x_n)}  {q_n}   \right\rfloor  \right\}     }
\end{equation}
where $\nu_p$ is the valuation corresponding to the prime $p$. We   call a point $\p \in \wP_\w^n (k)$ a \textbf{normalized point} if the weighted greatest common divisor of its coordinates is 1. 
 For any point $\p \in \wP_w^n (k)$, there exists its normalization given by  $ \q = \frac 1 {\wgcd (\p)} \star \p$.
Moreover, this normalization  is unique up to a multiplication by a $q$-root of unity, where $q=\gcd (q_0, \dots, q_n)$, see \cite{b-g-sh}. 
The \textbf{absolute weighted greatest common divisor} of    $\x=(x_0, \dots , x_n)$   is the largest 
real number $d$ such that %
$d^{q_i} \in \Z$    and $d^{q_i}\, \mid \, x_i$,  for all      $i=0, \dots n$. 
  We  denote it    by $\awgcd (x_0, \dots, x_n)$.

%
%Let $\w=(q_0, \dots , q_n)$ be a set of weights and $\wP^n(k)$ the weighted  projective space  over  a  field $k$.   
Let  $\p=[x_0, \dots , x_n] \in \wP^n(k)$.   Without any loss of generality we can assume that $\p$ is normalized.    The  \textbf{weighted multiplicative height} of $\p$   is  
\begin{equation}\label{def:height}
\wh_k( \p ) := \prod_{v \in M_k} \max   \left\{   \frac{}{}   |x_0|_v^{\frac {n_v} {q_0}} , \dots, |x_n|_v^{\frac {n_v} {q_n}} \right\}
\end{equation}
The height $\wh_k(\p)$ is well defined, in other words it does not depend on the choice of  coordinates of $\p$ and  $\wh_k(\p) \geq 1$; see \cite{b-g-sh}.
 Denote by $K = k (\awgcd (\p))$.  Then,  over $K$, the weighted greatest common divisor is the same as the absolute greatest common divisor,   $ \wgcd_K (\p) = \awgcd_K (\p)$.    Moreover, $[K:k] < \infty$ and we have the following. 

%****************************
 
\begin{prop}[\cite{b-g-sh}]\label{prop-5}
%Let $\w$, $K=\Q(\awgcd (\p) )$,  and $\p \in \wP^n(K)$, say $\p=[x_0 : x_1 :  \ldots  : x_n]$.      Then the following are true:
 If       $\p$ is normalized in $K$, then
\begin{equation}
\wh_K (\p)=  \wh_\infty (\p) = \max_{0 \leq i \leq n}\left\{\frac{}{}|x_i|^{{n_{\nu}}/q_i}_\infty \right\}.
\end{equation}
Moreover, if   $L/K$ is a finite extension,  then $\wh_L(\p)=  \wh_K (\p)^{[L:K]}$. 
\end{prop}

Using \cref{prop-5},   we can define the height on $\wP^n (\overline \Q)$. The height of a point on $\wP^n(\overline \Q)$ is called the \textbf{absolute (multiplicative) weighted height} and is the function 
\begin{equation}
\awh: \wP^n(\bar \Q)  \to [1, \infty),
\end{equation}
such that 
\begin{equation}
\awh(\p) =\wh_K(\p)^{1/[K:\Q]},
\end{equation} 
where  $\p \in \wP^n(K)$, for any $K$ which contains $\Q (\awgcd (\p))$. 
Moreover, for  $\p \in \wP^n(\overline \Q)$ and $\sigma \in G_{ \Q}$ we have $\wh (\p^\sigma) = \wh (\p)$. 
The \textbf{field of definition} of $\p$ is defined as 
$
\Q (\p) := \Q \left(       \left( \frac {x_0} {x_i} \right)^{\frac {q_0} {q}},  \ldots  , 1,  \ldots  ,   \left( \frac {x_n} {x_i} \right)^{\frac {q_n} {q}}  \right)
$
For any point $\p \in \wP_\w^n (\overline \Q)$, we have  $[ \Q (\p) : \Q] \leq  q \cdot [\Q(\phi (\p) ) : \Q ]$. 
The following result is  analogue to  Northcott's theorem for weighted projective spaces; see \cite{b-g-sh} for the proof. 

\begin{thm}[\cite{b-g-sh}] \label{thm_finite}
Let $c_0, d_0 \in \R$.
%   and $\wP_w^n(\overline \Q)$ the weighted projective space with weights  $\w = (q_0, \dots, q_n)$. 
Then the set 
\begin{equation}
\{\p \in \wP_w^n(\overline \Q): \wh_\Q (\p) \leq c_0 \text{ and } [\Q(\p):\Q] \leq d_0\}
\end{equation}
contains only finitely many points. 
\end{thm} 

% from here on the number field becomes K.  Change it to k

Hence,    $\{\p \in \wP_w^n( \overline \Q): \wh_{\bar \Q} (\p) \leq c_0 \} $ is a finite set for any constant $c_0$.    For any number field $k$, the set 
$ \{\p \in \wP_w^n(k):  \; \Q (\p) \subset k \; \text{ and } \;       \wh_k(\p) \leq c_0 \}$, 
is a finite set.    The next result is the analogue of Kronecker's theorem; see  \cite{b-g-sh}.

\begin{thm}[\cite{b-g-sh}]
Fix any $i$ with $x_{i} \neq 0$. Then $\wh (\p)= 1$ if  the ratio $x_j/\xi_{i}^{q_j}$, where $\xi_i$ is the $q_i$-th root of unity of $x_i$, is a root of unity or zero for every $0 \leq j \leq n$ and $j \neq i$.
\end{thm}

%*************** end of heights 
%\subsection{Curves with minimal invariants}

Now   we have the following two problems in terms of curves.  
%Let $\X$ be  a given superelliptic curve  with equation $z^m = f(x)$, $\deg f = d$, defined over  $\O_k$, and with corresponding moduli point $\p \in \wP_\w^n (k)$. 

\begin{prob}
 Given a curve $\X: z^m y^{d-m} = d(x, y)$ defined over $\O_k$,  determine     $\X^\prime$, $k$-isomorphic to $\X$, such that 
 defined over $\O_k$,   say $\X^\prime:   z^m y^{d-m} = g(x, y)$,  such that  $\p := \left[ \I (g) \right]  \in \wP_\w^n  (k)$  has minimal height over $k$. 
\end{prob}

\begin{prob}
Determine a twist $\Y$  of $\X$ such that $\p := \left[ \I (\Y) \right]  \in \wP_\w^n  (k)$  has minimal height over the algebraic closure $\bar k$.
\end{prob}

The above problems are equivalent of finding a model for the superelliptic curve such that the corresponding weighted moduli point has minimal possible weighted height or finding a twist with such property. 
% In the next section we will explain in more detail what this  means. 

%**************************
%% \newpage
%*****************
\section{Reduction of the moduli point}
The reduction of superelliptic curves consists of two steps.    
On the first step we perform the necessary coordinate changes so we have a minimal weighted moduli point, and the second step is to minimize  the coefficients of the equation of the curve.  Both of these steps are possible due to the fact that a superelliptic curve can be written as a curve with projective equation  $y^n z^{d-n}=f(x, z)$, such that $f$ is a binary form of degree  $\deg f =d$ with nonzero discriminant.  

The first step is based on the concept of weighted heights in weighted projective spaces as defined in \cite{b-g-sh}.   
%Applications of this method were used in \cite{mandili, b-guest, b-polak, rat-pts}.  
The second step is well known by work of Hermite for quadratics and extended by Julia for higher degree forms. For more recent work on this type of reduction see work of Cremona,  Stoll \cite{stoll-cremona} and Beshaj \cite{beshaj, beshaj-thesis, beshaj-2}.

%********************

%\newpage
\subsection{Minimal integral models of binary forms}\label{min-model}
%Integral binary forms with smallest invariants}
%
We say that a binary form $f(x, y )$ has a \textbf{integral minimal model} over $k$ if it is integral (i.e.  $f\in \O_k [x, y]$)   and  $\val_p (\I (f))$ is minimal for every prime $p \in \O_k$.

Let $f \in \O_k$ and $\x:= \I (f) \in \wP_\w^n (\O_k)$ its corresponding weighted moduli point. We define the \textbf{weighted valuation} of the tuple $\x=(x_0, \dots x_n)$ at the prime  $p\in \O_k$ as 
\begin{equation}
\val_p (\x) := \max \left\{   j   \, | \,  p^j \text{ divides }   x_i^{q_i}  \text{ for all }  i=0, \ldots n    \right\},
\end{equation}
Then we have the following. 
\begin{prop} \label{prop-6}
A binary form $f\in V_d$ is   a minimal model over $\O_k$ if for every prime $p \in \O_k$ such that $p \, | \, \wgcd (\I(f))$ the following holds 
\begin{equation} \val_p ( \I (f) ) < \frac d 2 q_i,  \quad \text{ for all } \; \; i=0, \dots , n.
\end{equation} 
 Moreover, for every integral binary form $f$ its minimal model exist.  
\end{prop}

\proof  Let $\x=\I(f)$.
From \cref{prop-2} we know that for any $M\in \GL_2(\O_k)$, $\I (f^M) = \left(  \det M \right)^{\frac d 2 } \I (f)$. 
   Hence, for every prime $p \in \O_k$ which divides $\wgcd (\x)$ we must "multiply"  $\x$ by the maximum exponent $j$ such that $\left( p^{\frac d 2 } \right)^j$ divides $\wgcd (\x)$.  

For a given binary form $f$ we pick $M= \begin{bmatrix}  \frac 1 \lambda  & 0 \\ 0 & 1  \end{bmatrix} $, where  
$\lambda$ is the weighted greatest common divisor of $\I(f)$ with respect the weights   $  \left(\left\lfloor  \frac {dq_0} {2} \right\rfloor,\ldots, \left\lfloor  \frac {dq_n} {2} \right\rfloor \right)$. 
The transformation $x \to \frac x \lambda$ gives a minimal model of  $f$ over $\O_k$. 
This completes the proof. 
\qed 

%Jaime
\begin{rem} If a prime $p \in \O_k$ divides  $\wgcd (\x)$ then $p^{q_i}$ divides $x_i$, so $ p^{q_i^{2}}$ divides   $x_i^{q_i}$.  Taking  $q_t=\min(q_0,\ldots,q_n)$, we have a lower bound  for the weighted valuation  of the 
point $\I (f)=(x_0,\ldots,x_n)$, that is  $\val_p ( \I (f) )\geq q_t^2$.
\end{rem}
Notice that it is possible to find a twist of $f$ with "smaller" invariants. In this case the new binary form is not in the same $\GL_2 (\O_k)$-orbit as $f$.  For example, the transformation 
$(x, y) \to \left( \frac 1 {\lambda^{\frac 2 d}} x,  \frac 1 {\lambda^{\frac 2 d}} y    \right)$,
will give us the form with smallest invariants, but not necessarily $k$-isomorphic to $f$. 

It is worth noting that for a binary form $f$ given in its  minimal model, the point $\I (f)$ is not necessarily normalized as in the sense of \cite{b-g-sh}.  

\begin{cor}
If $f(x, y) \in \O_k[x, y]$ is a binary form such that $\I(f)\in \wP_\w^n (k)$ is normalized over $k$, then $f$ is a minimal model over $\O_k$. 
\end{cor}

We see an example for binary sextics. 

\begin{exa}\label{exa-1}
Let be given the sextic 
\begin{equation} f(x, y) = 7776  x^6+31104  x^5 y +40176  x^4 y^2+25056  x^3 y^3+8382  x^2 y^4+1470  x y^5 +107 y^6\end{equation}
Notice that the polynomial  has content 1, so there is no obvious substitution here to simplify sextic. The moduli point is $\p=[J_2: J_4: J_6:J_{10}]$, where 
\begin{equation}
\begin{split}
J_2 & = 2^{15} \cdot 3^5, \; J_4  = - 2^{12}    \cdot 3^9 \cdot   101 \cdot  233,  \;   J_6  = 2^{16}  \cdot 3^{13} \cdot 29 \cdot 37 \cdot 8837, \\
J_{10} & = 2^{26} \cdot 3^{21}  \cdot 11  \cdot 23 \cdot  547 \cdot 1445831\\
\end{split}
\end{equation}
 Recall that the transformation $(x, y) \to \left( \frac 1 p x , y\right)$ will change the representation of the point $\p$ via 
\begin{equation} \frac 1 {p^3} \star [J_2: J_4: J_6:J_{10}]= \left[ \frac 1 {p^6}   J_2: \frac 1 {p^{12}} J_4: \frac 1 {p^{18}} J_6:   \frac 1 {p^{30}}J_{10} \right]
\end{equation}
So we are looking for prime factors  $p$ such that $p^6|J_2$, $p^{12} | J_4$, $p^{18}| J_6$, and $p^{30} | J_{10}$.  Such candidates for $p$ have to be divisors of 
 $\wgcd (\p) = 2^2 \cdot 3^2$. 

Obviously  neither $p=2$ or $p=3$ will work.  Thus, $f(x, y)$ is in its minimal model over $\O_k$. 
%
\iffalse
However, the transformation 
%
\begin{equation} 
(x, y) \to \left( \frac 1 {6^{\frac 1 3}} x, \frac 1 {6^{\frac 1 3}} y  \right).
\end{equation}
%
will give us 
%
\begin{equation}
g(x, y) = 216\,{x}^{6}+864\,{x}^{5}y+1116\,{x}^{4}{y}^{2}+696\,{x}^{3}{y}^{3}+{
\frac {1397\,{x}^{2}{y}^{4}}{6}}+{\frac {245\,x{y}^{5}}{6}}+{\frac {
107\,{y}^{6}}{36}}
\end{equation}
%
with set of invariants 
%
\begin{equation} 
I(g)= [ 2^{11} \cdot 3: - 2^4 \cdot 3 \cdot  101 \cdot  233: 2^4 \cdot 3 \cdot 29 \cdot 37 \cdot 8837: 2^6 \cdot 3  \cdot 11  \cdot 23 \cdot  547 \cdot 1445831], 
\end{equation}
%
which is obviously  normalized in $\wP_\w^3 (\Q)$ since $\wgcd (\p^\prime) =1$.
%
\fi
%
\qed
\end{exa}

\begin{cor} 
The transformation of $f(x, y)$ by the matrix 
\begin{equation} 
M= \begin{bmatrix}
\varepsilon_d \frac 1 {\left( \wgcd (I(f)) \right)^{\frac 2 d}} & 0 \\ 0 & \varepsilon_d \frac 1 {\left( \wgcd (I(f)) \right)^{\frac 2 d}}
\end{bmatrix}
\end{equation}
where $\varepsilon_d$ is a $d$-primitive root of unity, will always give a minimal set of invariants.
\end{cor}

%**********************
\subsection{Reduction a):  Moduli points with minimal weighted height}\label{red-a}
Let $\X$ be as in \cref{w-eq-super}  and $\p = [\I(f)]\in \wP_\w^n (k)$.  
Let us assume that for a prime $p\in \O_k$, we have $\nu_p \left( \wgcd (\p) \right) = \alpha$.  If we use the transformation 
$ x\to \frac x {p^\beta} x, \quad \text{ for }  \quad \beta \leq \alpha$, 
 then from \cref{prop-2}   the set of invariants will become 
$ \frac 1 {p^{ \frac d 2 \beta }}\star \I(f)$.
To ensure that the moduli point $\p$ is still with integer coefficients we must pick $\beta$ such that $p^{\frac {\beta d} 2}$ divides  $p^{\nu_p (x_i)}$ for $i= 0, \dots , n$.   Hence, we must pick $\beta$ as the maximum integer such that $\beta \leq \frac 2 d \nu_p (x_i)$, for all $i=0, \dots , n$. This is the same $\beta$ as in \cref{prop-6}. 
The transformation 
$ (x, y) \to \left( \frac x {p^\beta} , y \right)$, 
has corresponding matrix $M=\begin{bmatrix} \frac 1 {p^\beta} & 0 \\ 0 & 1 \end{bmatrix}$ with $\det M = \frac 1 {p^\beta}$.  Hence, from \cref{prop-2} the moduli point $\p$ changes as 
$  \p \to \left( \frac 1 {p^\beta} \right)^{d/2} \star \p $,  
which is still an integer tuple.  We do this for all primes $p$ dividing $\wgcd (\p)$.     Notice that the new point is not necessarily normalized in $\wP_\w^n (k)$ since $\beta $ is not necessarily equal to $\alpha$.  
This motivates the following definition. 

\begin{defi}
Let $\X$ be a superelliptic curve defined over an integer ring $\O_k$ and $\p \in \wP_\w^n (\O_k)$ its corresponding weighted moduli point. We say that $\X$ has a \textbf{minimal model} over $\O_k$ if for every prime $p \in \O_k$ the \textbf{valuation of the tuple} at $p$
\begin{equation}
\val_p (\p) := \max \left\{  \nu_p (x_i)   \text{ for all }  i=0, \ldots n    \right\},
\end{equation}
is minimal, where $\nu_p (x_i)$ is the valuation of $x_i$ at the prime $p$. 
\end{defi}
%
%***************
\begin{thm}\label{thm-1}  Minimal models of superelliptic curves exist.  
An  equation $\X :$ $ z^m y^{d-m}= f(x, y)$ is a minimal model over $\O_k$,   if   for every prime $p \in \O_k$ which divides  $ p \, | \,  \wgcd \left( \I(f) \right)$,  the valuation $\val_p $  of $\I(f)$ at $p$ satisfies  
\begin{equation}\label{val}
\val_p (\I(f))  <   \frac d 2  \,  q_i,
\end{equation}
for all $i=0, \dots , n$.  Moreover,  then for  $\lambda=\wgcd (\I (f)) $ with respect  the weights
$(\left\lfloor  \frac {dq_0} {2} \right\rfloor,\ldots, \left\lfloor  \frac {dq_n} {2} \right\rfloor) $
the transformation 
$ (x, y, z) \to \left(      \frac x \lambda, y, \lambda^{\frac d m} z \right) $
 gives a minimal model of  $\X$ over $\O_k$. If $m|d$ then this isomorphism is defined over $k$. 
\end{thm}

\proof   Let $\X$ be a superelliptic curve given by \cref{w-eq-super} over $\O_k$ and   $\p = \I(f) \in \wP_\w^n (\O_k)$ with weights $\w=(q_0, \dots , q_n)$.  Then $\p \in \wP_\w^n (\O_k)$ and from \cref{prop-6} exists $M\in \GL_2 (\O_k)$ such that $M= \begin{bmatrix}  \frac 1 \lambda  & 0 \\ 0 & 1  \end{bmatrix} $ and $\lambda$ as in the theorem's hypothesis.  By \cref{prop-6} we have that \cref{val} holds. 

Let us see how the equation of the curve $\X$ changes when we apply the transformation by $M$. We have 
%For $\lambda $ as in \cref{lambda} we change the coordinates as $(x, y, z) \to (\frac x \lambda, y, z)$.   The equation of the curve becomes
%
\begin{equation} 
z^m y^{d-m} = f  \left(\frac x \lambda , y \right) =  a_d \frac {x^d} {\lambda^d} + a_{d-1} \frac {x^{d-1}} {\lambda^{d-1}} y + \cdots + a_1 \frac x \lambda y^{d-1} + a_0 y^d
\end{equation}
Hence, 
\begin{equation}\label{eq-2}
\X^\prime : \;  \lambda^d  z^m y^{d-m} = a_d  x^d + \lambda a_{d-1}  x^{d-1}  y + \cdots + \lambda^{d-1} a_1 x   y^{d-1} + \lambda^d a_0  y^d
\end{equation}
This equation has coefficients in $\O_k$.  Its weighted moduli point is 
$ \I (f^M) = \frac 1 {\lambda^{\frac d 2}}\star \I(f)$, 
which satisfies \cref{val}.   It is a twist of the curve $\X$ since $\lambda^d$ is not necessary a $m$-th power in $\O_k$. The isomorphism of the curves over the field  $k \left(\lambda^{\frac d m}\right)$ is given by   $(x, y, z) \to \left( \frac x \lambda, y, \lambda^{\frac d m} z   \right) $.
If $m|d$ then this isomorphism is defined over $k$  and $\X^\prime$ has equation 
\begin{equation} 
\X^\prime : \;   z^m  y^{d-m} =  a_d  x^d + \lambda a_{d-1}  x^{d-1}  y + \cdots + \lambda^{d-1} a_1 x   y^{d-1} + \lambda^d a_0  y^d
\end{equation}
\qed

Then we have the following.

\begin{cor}
There exists a  curve $\X^\prime$ given in \cref{eq-2}  isomorphic to $\X$ over the field $K:=k \left( \wgcd( \p)^{\frac d m}\right)$ with minimal    invariants.    Moreover, if $m  \mid  d$ then $\X$ and $\X^\prime$ are $k$-isomorphic.
\end{cor}

A simple observation from the above is that in the case of hyperelliptic curves we have $m=2$ and $d=2g+2$.  Hence,  the curves $\X$ and $\X^\prime$ would always be isomorphic over $k$.  So we have the following. 

\begin{cor} 
Given a hyperelliptic curve  defined over a ring of integers $\O_k$. There exists a  curve $\X^\prime$  $k$-isomorphic to $\X$   with minimal    invariants.  
\end{cor}

%*****************
% 
\cref{thm-1} above provides an algorithm which is described next.
Given a superelliptic curve $\X$  defined over $\O_k$,  we denote the corresponding point in the weighted moduli space by $\p = \left[ x_0 : \dots : x_n    \right]$.  \\

\noindent  \textbf{Algorithm:} Computing an equation of the curve with  minimal moduli point. \\  

\noindent \textsc{Input:} A  curve $\X: z^m y^{d-m} = f(x, y)$,   $\deg (f)=d$ and  $f \in \O_k[x, y]$. \\
\noindent \textsc{Output:} A curve   $\Y : \lambda^d z^m = g(x, y) $,  defined over $\O_k$ and $k$-isomorphic to $\X$ such that    $\val_p \, \I(g)$ is minimal  for each prime  $p \in \O_k$. \\

\medskip

  \textsc{Step 1:}  Compute the generating set $\I:=[I_{q_0}, \ldots , I_{q_n}]$  for  $\cR_n$. 
  
   \textsc{Step 2:}  Compute the  moduli point $\p\in \wP_w^n (k)$  for $\X$ by evaluating $\I(f)$. 
   
  \textsc{Step 3:}  Computing $\lambda= \wgcd ( \I (f) ) \in \O_k $ with respect the weights $\left   (\left\lfloor  \frac {dq_0} {2} \right\rfloor,\ldots, \left\lfloor  \frac {dq_n} {2} \right\rfloor   \right) $

   \textsc{Step 4:}  Compute $f^M$,  where   $M : =\begin{bmatrix} \frac 1 \lambda & 0 \\ 0 & 1 \end{bmatrix}$. We have 
$ g(x) := \lambda^d \cdot f \left( \frac x \lambda    \right)$.

   \textsc{Step 5:} Return the curve $\X^\prime$  with   equation   $\lambda^d z^m y^{d-m} = g(x, y)$.  

\medskip

\noindent Let us illustrate  for $g=1$, 2.  \\

\subsubsection{Elliptic curves}\label{section-4-3}
Technically elliptic curves are not superelliptic curves, but the method will work the same.   Let $E$ be an elliptic curve with Weierstrass equation as in Birch/Swinnerton-Dyer \cite{birch-s-d-1}
\begin{equation}\label{ell-curve}
E : \; \; z^2 y^2 = f(x, y)= a x^4 + b x^3 y + ex^2 y^2 + c x y^3 + d y^4.
\end{equation}
Invariants of the binary quartic  $f(x, y)$ are  
\begin{equation}
I_2=12ae - 3bd+c^2, \; \; I_3= 72 ace + 9 bcd - 27ad^2-27eb^2-2c^3
\end{equation}
The corresponding weighted moduli space is $\wP_{(2, 3)} (k)$.   Isomorphism classes of elliptic curves over $k$ correspond to points in $\wP_{(2, 3)} (k)$. 

\begin{cor}
The equation in \cref{ell-curve} is a  minimal models if for every prime $p \in \Z$, $p\neq 2, 3$ which divides $p  | \wgcd( I_2, I_3)$, the valuation $\val_p (I_2, I_3)$ satisfies   $\val_p (I_2, I_3) < 2 q_i$, 
for $q_i=2, 3$. 
\end{cor}

Hence, we can reduce any prime $p\neq 2, 3$ such that $p^\a | I_2$ and $p^\b | I_3$ when $\a \geq 4$ and $\b \geq 6$. 
The above result was proved in \cite{birch-s-d-1} using a case by case analysis; see \cite[Lem. 3]{birch-s-d-1}.  Lemma 4 and Lemma 5 in \cite{birch-s-d-1} describe the cases when $p=2$ and $p=3$ respectively. 
Then we have the following; see \cite[Theorem 1]{birch-s-d-1}. 

\begin{thm}[Birch, Swinnerton-Dyer]
If $E$ is an elliptic curve,  and $U$ is a non-trivial 2-covering of $E$ then $U$ can be represented by a curve $z^2 y^2=g(x, y)$, where $g(x, y)$ is reduced. The reduced $g(x, y)$ of a given $E$ are finite in number and computable. 
\end{thm}

\subsubsection{Genus 2 curves}
Let $\X$ be a genus 2 curve with equation 
$ z^2 y^4 = f(x, y)$ as in \cref{exa-1}. 
By applying the transformation 
$ (x, y, z) \to \left( \frac x 6, y, 6^3\cdot z \right)$ 
 we get the equation 
\begin{equation}\label{eq-exa} 
z^2= x^6+24  x^5+186  x^4+696  x^3+1397  x^2+1470  x+642.
\end{equation}
Computing the moduli point of this curve we get 
\begin{equation} 
\p= [ 2^{11} \cdot 3: - 2^4 \cdot 3 \cdot  101 \cdot  233: 2^4 \cdot 3 \cdot 29 \cdot 37 \cdot 8837: 2^6 \cdot 3  \cdot 11  \cdot 23 \cdot  547 \cdot 1445831], 
\end{equation}
which is obviously  normalized in $\wP_\w^3 (\Q)$ since $\wgcd (\p) =1$. Hence, the \cref{eq-exa} is a minimal model.

%*******************
\subsection{Minimal integral twists}
Let $k$ be a field of characteristic zero and $\Gal (\bar k/k)$ the Galois group of $\bar k/k$. 
Let  $\X$ a genus $g\geq 2$ smooth, projective algebraic curve defined over $k$. We denote by $\Aut (\X)$ the automorphism group of $\X$ over the algebraic closure $\bar k$. By $\Aut_k (\X)$ is denoted the subgroup of automorphisms of $\X$ defined over $k$.  
A \textbf{twist} of $\X$ over k is a smooth projective curve $\X^\prime$ defined over $k$ which is isomorphic to $\X$ over $\bar k$.  We will identify  two twists which are isomorphic over $k$. The set of all twists of $\X$, modulo $k$-isomorphism, is denoted by $\twist (\X/k)$.

Let $\X$ and $\X^\prime$ be twists of each other over $k$. Hence, there is an isomorphism 
$ \phi: \X^\prime \to \X  $
defined over $\bar k$. 
For any $\sigma \in  \Gal (\bar k / k)$ there exists the induced map $\phi^\sigma : \X^\prime \to \X$. 
To measure the failure of $\phi$ being defined over $k$ one considers the map
$
\xi : \; \;  \Gal (\bar k / k)   \to \Aut (\X) 
$
such that $\xi (\sigma) =  \phi^\sigma \phi^{-1} $ for any $\sigma \in \Gal (\bar k / k)$.
The following is the main result on twists; see \cite[Th.~2.2. pg. 285]{silv-1} % and the proof can be found in \cite[Thm.~2.2]{silv-1}. pg 285
\begin{prop}  The following are true:
\begin{enumerate}[\upshape(i), nolistsep]
\item $\xi$ is a 1-cocycle
\item  The cohomology class $\{ \xi \}$ is determined by the $k$-isomorphism class of $\X^\prime$ and is independent of $\phi$.  Hence, there is a natural map
\begin{equation} \theta: \twist (\X /k) \to \homo^1 \left(\Gal (\bar k / k ), \Aut (\X) \right)  \end{equation}
\item  The map  $\theta$ is a bijection.
\end{enumerate}
\end{prop} 

As noted by Silverman in \cite[Remark 2.3, pg. 285]{silv-1}, $\homo^1 \left(\Gal (\bar k / k ), \Aut (\X) \right) $ is not necessarily a group since $\Aut (\X)$ is not necessarily Abelian. For nonabelian Galois cohomology we refer to \cite{serre}.

%The focus of this paper is on superelliptic curves.  Next we see how the general theory of twists is applied to such curves. 

Fix an integer $n \geq 2$.  Let $k$ be a number field which contains  all primitive $n$-th roots of unity $\xi_n$  and $\O_k$ its ring of integers. Consider    a smooth superelliptic curve  $\X$ defined over $\O_k$ with equation 
$ y^n = f(x) $
where $f$ is a separable   polynomial in $k$.   The multiplicative group of $n$-th roots of unity will be denoted by $\mu_n$.
Then, $\mu_n$ is embedded in $\Aut (\X)$ in the obvious way. 
Hence, $\Aut (\X)$ is an extension of the group $\mu_n$.  The  list of isomorphism classes of such groups is determined; see  \cite{aut}.  Thus, determining the set $\twist (\X/k)$ is equivalent to determining $\homo^1 \left(\Gal (\bar k / k ), \Aut (\X) \right) $ for each possible group $\Aut (\X)$. 

%****************************************************************
\begin{lem}
Let $\X$ be as in \cref{w-eq-super}.
Any twist $\X^\prime$ of $\X$ has equation 
$ \lambda  y^n = f(x)$, 
for some $\l\in k$ such that $\lambda^{1/n} \in k$.  Moreover, the isomorphism $\phi : \X  \longrightarrow  \X^\prime$ is given by
$  (x, y)  \longrightarrow (x, \sqrt[n]{\l} \, y) $.
\end{lem}

\proof
Two curves  $\X$ and  $\X^\prime$,   with equations
$y^n = f(x)$  and $ v^n = g(u)$,      where $\deg f = \deg g =s$,   are isomorphic   over $\bar k$, if and only if 
\begin{equation}
x= \frac {au+b} {cu+d}, \; y= \frac {\l \, v } {(cu+d)^{s/n}}, \; \; \text{where } \begin{bmatrix} a & b \\ c & d \end{bmatrix} \in \GL_2 (k), \; \l \in k^\ast 
\end{equation}
Let $f(x) $ be given by  $ f(x) =\sum_{i=0}^sa_i x^i$ % a_0 + a_1 x + \cdots + a_{s-1} x^{s-1} + a_s x^s, \end{equation}
then 
\begin{align}
 f \left( \frac {au+b} {cu+d} \right) & = a_0 + a_1 \left( \frac {au+b} {cu+d} \right) + \cdots   + a_s \left( \frac {au+b} {cu+d} \right)^s  =  \frac 1 {(cu+d)^s } \, g(u), 
\end{align}
where $g(u)$ is of degree $s$ in $u$.   hence, the equation of the curve becomes   $(cu+d)^s  \, y^n = g(u) $.
Replacing $y$ as above we get    $   \l^n \, v^n   = g(u)$. 
\qed 

%***************
\begin{thm}\label{thm-2}  
Minimal twists of minimal models of superelliptic curves exist.  
An  equation $\X :$ $ z^m y^{d-m}= f(x, y)$ is a minimal twist  over $\O_k$,   if   for every prime $p \in \O_k$ which divides  $ p \, | \,  \awgcd \left( \I(f) \right)$,  the valuation $\val_p $  of $\I(f)$ at $p$ satisfies  
\begin{equation}\label{val-twist}
\val_p (\I(f))  <   \frac d 2  \,  q_i,
\end{equation}
for all $i=0, \dots , n$.  Moreover,  then for  $\lambda=\wgcd (\I (f)) $ with respect  the weights
$(\left\lfloor  \frac {dq_0} {2} \right\rfloor,\ldots, \left\lfloor  \frac {dq_n} {2} \right\rfloor) $%
the transformation 
$ (x, y, z) \to \left(      \frac x \lambda, y, \lambda^{\frac d m} z \right) $
 gives a minimal model of  $\X$ over $\O_k$. If $m|d$ then this isomorphism is defined over $k$. 
\end{thm}

\proof   Let $\X$ be a superelliptic curve given by \cref{w-eq-super} over $\O_k$ and   $\p = \I(f) \in \wP_\w^n (\O_k)$ with weights $\w=(q_0, \dots , q_n)$.  Then $\p \in \wP_\w^n (\O_k)$ and from \cref{prop-6} exists $M\in \GL_2 (\O_k)$ such that $M= \begin{bmatrix}  \frac 1 \lambda  & 0 \\ 0 & 1  \end{bmatrix} $ and $\lambda$ as in the theorem's hypothesis.  
%By \cref{prop-6} we have that \cref{val} holds. 

Let us see how the equation of the curve $\X$ changes when we apply the transformation by $M$. We have 
%For $\lambda $ as in \cref{lambda} we change the coordinates as $(x, y, z) \to (\frac x \lambda, y, z)$.   The equation of the curve becomes
%
\begin{equation} 
z^m y^{d-m} = f  \left(\frac x \lambda , y \right) =  a_d \frac {x^d} {\lambda^d} + a_{d-1} \frac {x^{d-1}} {\lambda^{d-1}} y + \cdots + a_1 \frac x \lambda y^{d-1} + a_0 y^d
\end{equation}
Hence, 
\begin{equation}\label{eq-2-twist}
\X^\prime : \;  \lambda^d  z^m y^{d-m} = a_d  x^d + \lambda a_{d-1}  x^{d-1}  y + \cdots + \lambda^{d-1} a_1 x   y^{d-1} + \lambda^d a_0  y^d
\end{equation}
This equation has coefficients in $\O_k$.  Its weighted moduli point is 
\begin{equation} \I (f^M) = \frac 1 {\lambda^{\frac d 2}}\star \I(f), \end{equation}
which satisfies \cref{val-twist}.   It is a twist of the curve $\X$ since $\lambda^d$ is not necessary a $m$-th power in $\O_k$. The isomorphism of the curves over the field  $k \left(\lambda^{\frac d m}\right)$ is given by 
$ (x, y, z) \to \left( \frac x \lambda, y, \lambda^{\frac d m} z   \right) $.
If $m|d$ then this isomorphism is defined over $k$  and $\X^\prime$ has equation 
\begin{equation} 
\X^\prime : \;   z^m  y^{d-m} =  a_d  x^d + \lambda a_{d-1}  x^{d-1}  y + \cdots + \lambda^{d-1} a_1 x   y^{d-1} + \lambda^d a_0  y^d
\end{equation}
This completes the proof.
\qed

Further work on the reduction of moduli points of higher degree binary forms is intended in \cite{curri}.

%***************************************************
%\newpage
\section{Reduction of coefficients of binary forms}\label{section-4-2}
Next we will focus on the reduction of coefficients of a binary form.  This is an old problem starting with Julia's thesis in \cite{julia} and continued with more recent papers \cite{beshaj-thesis}, \cite{beshaj}, \cite{heights}, \cite{beshaj-1}, \cite{beshaj-2}. 
Along similar lines a new reduction of binary forms via the hyperbolic centroid
 was introduced in \cite{e-sh}, which seems to have different results from the approach in previous papers (cf. \cref{hyp-centroid}). 

%************************************************
%\subsection{Heights of polynomials}\label{heights_pol}     

A non-homogenous polynomial with $n$ variables will be denoted as  
\begin{equation}
f(x_1, \dots, x_n) = \sum_{\substack{i=(i_1, \dots, i_n) \in I} } a_i x_1^{i_1} \cdots x_n^{i_n},
\end{equation} 
where all $a_i \in k$,  $I \subset \left(  \Z^{\geq 0}\right)^n$, and $I$ is finite. Let $\deg f$ denote the total degree of $f$. We will use lexicographic ordering to order the terms in a given polynomial, and $x_1 > x_2 > \dots >x_n$.   The \textbf{(affine) multiplicative height of $f$} is defined as  
\begin{equation}
\h^{\mathbb A}_k(f)= \prod_{v \in M_k} \max\left\{\frac{}{}1, |f|_v^{n_v}\right\}, \; \; 
\end{equation}
where     $| f |_v := \max_j\left\{\frac{}{} |a_j|_v \right\}$   is called the \textbf{Gauss norm} for any absolute value $v$ and     
 $n_v$ is the \textbf{local degree at $v$} given by $n_v= [ k_v : \Q_v]$ for $k_v$ and $\Q_v$ are the completions with respect to $v$; see \cref{heights}. 
Hence, the affine height of a polynomial is defined to be the height of its coefficients taken as affine coordinates.  

The \textbf{(projective) multiplicative height of a polynomial} is the height of its coefficients taken as coordinates in the projective space. Thus,   
\begin{equation}
\h_k(f)= \prod_{v \in M_k}  | f |_v^{n_v}
\end{equation}
%   define n_v
The   \textbf{(projective) absolute multiplicative  height} is defined as 
\begin{equation}
 \h: \P^n(\Q)  \to [1, \infty),
 \end{equation}
 such that   $ \h(f)=\h_k(f)^{1/[k:\Q]}$.
%     and in the same way $h(f)$, $\h^{\mathbb A}(f)$, $h^{\mathbb A}(f)$. 
The following is \cite[Thm.~2, Prop.~1]{heights}.
\begin{lem}\label{pol_finite} The following hold true:
\begin{enumerate}[\upshape(i), nolistsep]
\item Let   $f \in k [x, y]$. Then  there are only finitely many polynomials $ g  \in k [x, y]$ such that $\h_k(g) \leq \h_k(f)$.  
\item Let $f(x_0, \dots, x_n)$ and $g(y_0, \dots, y_n)$ be polynomials in different variables. Then,   $\h(f \cdot g)= \h(f) \cdot \h(g)$.
\end{enumerate}
\end{lem}
The following lemma is true for the product of a finite number of polynomials. 
\begin{lem}[Gauss's lemma] 
Let $k$ be a number field and $f, g \in k[x_1, \dots, x_n]$.  If $v$ is not Archimedean, then $|f g|_v= |f|_v  |g|_v$. 
\end{lem}
The proof can be found in \cite[pg. 22]{bombieri}.  Gauss's lemma applies to all non-Archimedean absolute values but the Archimedean case is more complicated. 
The following   gives a bound for the homogenous polynomial evaluated at a point; see \cite[Lem.~15]{heights}.  

\begin{lem}\label{hom-at-point}  Let $k$ be a number field, $f\in k[x_0, \dots, x_n]$ a homogenous polynomial of degree $d$, and $\a=(\a_0, \dots, \a_n)\in \overline k^{n+1}$. Then  the following hold:
\begin{enumerate}[\upshape(i), nolistsep]
\item    $|f(\a)|_v \leq |c(d, n)|_v  \cdot \max_j \left \{\frac{}{}|\a_j|_v\right\}^d \cdot  |f |_v$,  where $| c(d, n) |_v$ is $\binom{n+d}{d}$ if $v$ is non-Archimedean and 1 otherwise.
\item  $ \h(f(\a)) \leq  c_0 \cdot  \h(\a)^d \cdot \h(f).$   
\end{enumerate}
\end{lem}
Then we have the  immediate corollary; see \cite[Cor.~1]{heights}.

\begin{cor}
Let $k$ be a number field, $f \in k[x,z]$ a homogenous polynomial of degree $d$ as $f(x, z) = \sum_{i=1}^d a_ix^iz^{d-i}$  and  $\a=(\a_0, \a_1) \in \overline k^{2}$. Then,
\begin{equation}\h(f(\a)) \leq  \min \left \{  d+1, 2^{d+1} \right\} \cdot  \h(\a)^d \cdot \h(f).
\end{equation}
\end{cor}
We will use Lemma~\ref{hom-at-point} to bound the height of the   invariants on $V_d$; see \cite[Thm.~4]{heights} for the proof.
\begin{thm}
Let $M =\left[ a_{i, j} \right]  \in \GL_2 (k)$,  $f\in k[x_1, \ldots , x_n]$,  $f \in V_d (k)$,  and $\h(f)$   the absolute height of $f$.  Then,
\begin{equation} 
\h( f^M ) \leq 2^n \cdot (n+1) \cdot \h(M)^n  \cdot   \h(f),
\end{equation}
where $\h(M)= \max \{ a_{i, j} \}$.
\end{thm}
Let $k$ be an algebraic number field, and $f (x, y)$ and  $\bar f(x, y) := f(u x + w , y)$.  Then from \cite[Thm.~5]{heights} we have:
\begin{thm}
% The following are true:
\begin{enumerate}[\upshape(i), nolistsep]
\item  For any valuation $v\in M_k$ we have 
\begin{equation} | \bar f |_v \leq 2_v^{d}  \, \cdot \, c(d)_v  \, \cdot \, |u|_v^d \, \cdot \, |w|_v^d \, \cdot \, \max_{0 \leq i \leq d}  \left\{ \frac{}{} |b_i|_v\right\} \end{equation}
\item The height of the form is bounded as follows
\begin{equation} \h( \bar f) \leq (d+1)  \, \cdot \,    2^{d} \, \cdot \, u^{d } \, \cdot \, w^{d }    \, \cdot \, \h(f) \end{equation}
\end{enumerate}
\end{thm}

Denote by $\Orb (f)$ the $\GL_2(k)$-orbit of $f$  in $V_d$ and $\h(f)$ its height.  Note that there are only finitely many $f^\prime \in \Orb(f)$ such that $\h(f^\prime) \leq \h(f)$.  Define the \emph{minimum height}  of the binary form $f(x, y)$ as follows
\begin{equation} 
\tilde \h(f) := \min \left\{ \frac{}{} \h(f^\prime) | f^\prime \in \Orb(f), \, \h(f^\prime) \leq \h(f) \right\}.
\end{equation}
We naturally have the following problem:
\begin{prob}
For every $f$ let $f^\prime$ be the binary form such that $f^\prime \in \Orb(f)$ and $\tilde \h(f) = \h(f^\prime)$. Determine a matrix $M \in \GL_2(k)$ such that $f^\prime = f^M$. 
\end{prob}
This can be fully solved for quadratics (cf. cref{binary-quad}).  There is a connection between the height of the moduli point $\p_f$  (considered as a projective point in $\P^n(k)$) and $\tilde \h(f) $ as described in \cite[Thm.~6]{heights}. We have 
\begin{equation}
\h (\p) \leq c \cdot \tilde \h(f), 
\end{equation}
for some constant $c$.  For binary sextics this constant was computed in \cite{heights} as $c=2^{28} \cdot 3^9 \cdot 5^5 \cdot 7 \cdot 11 \cdot 13 \cdot 17 \cdot 43$; see \cite[Lem.~20]{heights}.

%*******************************
\subsection{Quadratic forms} \label{binary-quad}
The case of quadratic forms is well known and goes back to Lagrange, Gauss \cite{gauss, gauss-2}, Hermite and many others.  
 A \textbf{quadratic form over $\R$} is a function $\q: \R^n \to \R$ that has the form $\q(\x)= \x^TA\x$  where $A$ is a symmetric $n \times n$ matrix called the \textbf{matrix of the quadratic form}.    Let $f$,  $g$ be quadratic forms and $A_f$,  $A_g$ their corresponding matrices. Then, $f \sim g$ if and only if $A_f$ is similar to   $A_g$. 

Let $\q(x, y)= ax^2+ bxy + cy^2$ be a binary quadratic in $\R[x, y]$. We will use the following  notation  to represent the binary quadratic, $\q(x, y) = [a, b, c]$. The \textbf{discriminant} of $\q$ is $D= b^2-4ac$ and  $\q(x, y) $ is positive definite if $a>0$ and $D <0$. Denote the \emph{set of positive definite binary quadratics} with $V_2^+ (\R)$.    Let $\SL_2(\R)$ acts   as usual on the set of positive definite binary quadratic forms $V_2^+ (\R)$. 
%

%***************************
%\subsubsection{The zero map}
%
Consider  the following map $\z: V_2^+(\R) \to \H_2$,     which is called the \textbf{zero map}:
\begin{equation}\label{zero-map-real}  
[a, b, c]  \to \z(\q)  = \frac{-b + \sqrt \D}{2a}.
\end{equation}
%
%where $\re(\z(\q)) =-\frac {b}{2a}$, and $\im (\z(\q)) =\frac{ \sqrt{|\D|}}{2a}$.   
It is a bijection since given $z=x+iy$, we can find $a, b, c$ such that $\q(x, y)$ is positive definite  given as   $ [1, -2x, x^2+y^2]$. 
 %\begin{lem} 
The map $\z$ gives us   a one to one correspondence between positive definite quadratic forms and points in $\H_2$. 
%\end{lem} 
Let $\Gamma$ be the modular group acting on $\H_2$ and on $V_2^+$ as described above.  
Then, from \cite{beshaj-thesis}, 
%\begin{lem}\label{equivariant-real}
the zero  map $ \z$ is a $\Gamma$-equivariant  map (i.e.  $\z  (\q^M)=M^{-1}   \z(\q)$, for any $M\in \Gamma$).
%\end{lem}

  Define $\q \in V_2^+ (\R)$ to be \textbf{reduced} if $\z(\q) \in \F$.  
%
%\begin{thm}\label{red-quad} 
The following are facts are well known.  We suggest \cite{beshaj-thesis} among many other references. 

\begin{enumerate}[\upshape(i), nolistsep]
\item   A  quadratic form $\q \in V_2^+(\R)$ is  reduced if and only if $|b| \leq a \leq c$. 
%\begin{thm}\label{bound-b}
\item Let $\q$ be  reduced  with   discriminant $\D=-D$. Then, $b \leq \sqrt{ D/ 3}$.
\item The number of reduced forms of a fixed discriminant  $\D=-D$  is finite.
%\end{thm}
%\begin{thm}\label{discrim-red-finite}
\item   Every  $\q \in V_2^+ (\R)$   is equivalent to a reduced form of the same discriminant. 
\end{enumerate}
%\end{thm}
Then we have the following; see \cite{beshaj-thesis, beshaj-2} among other places.
\begin{thm} \label{red-min-reals}
Let  $f(x, y) = ax^2+bxy+cy^2$ be reduced.  Then,  $\h( f) =c$.   Moreover,  $f$ obtains it absolute minimal height  in its $\Gamma$-orbit. 
\end{thm}
\begin{rem}
Lagrange proved that for every value $\D$, there are only finitely many classes of binary quadratic forms with discriminant $\D$. Their number is the \textbf{class number} of discriminant $\D$. He described an algorithm, called \emph{reduction}, for constructing a canonical representative in each class, the reduced form, whose coefficients are the smallest in a suitable sense.
Gauss \cite{gauss, gauss-2} gave a better reduction algorithm in Disquisitiones Arithmeticae.   For a more modern treatment for reduction of quadratics see also \cite{zagier}.
%, which ever since has been the reduction algorithm most commonly given in textbooks. In 1981, Zagier published an alternative reduction algorithm which has found several uses as an alternative to Gauss's.[5]
\end{rem}

%************* chap 4
\subsection{Reduction  of  higher degree binary forms in their $\G$-orbit}\label{red-bin-forms}
Let $f\in V_d (k)$  and $\G_{\O_k} := \SL_2 (\O_k) / \{ \pm I \} $.  To every binary form $f$ it is associated a positive definite quadratic $\J_f$ called the \textit{Julia quadratic}. In \cite{julia} is proved that this  is a covariant of the degree $d$ binary forms and we can develop reduction theory using this quadratic.    A degree $d$ binary form is called \textit{reduced} when $\z(\J_f)$  is in the fundamental domain of the action of the modular group $\G$ on $\H_2$. 
 
Let $f\in V_d (\R)$ as in \cref{w-eq-super}. We consider $f$ as a polynomial in a single variable $f(x,1)$.  Let     the real roots of $f(x,y)$ be $\a_i$,   for $1 \leq i \leq r$ and the pair of complex roots $\b_j$, $\bar \b_j$ for $1 \leq j \leq s$, where $r+2s =d$. 
Then 
\begin{equation}\label{deg-n-real-factored}
f(x,1) = \prod_{i=1}^r (x-\a_i) \cdot\prod_{i=1}^s (x-\b_i)(x-\bar \b_i).
\end{equation}
The ordered pair $(r, s)$ of numbers $r$ and $s$ is called the \textbf{signature} of  $f$.   We associate to $f$  two   quadratic forms, which we are writing as polynomials, namely  $T_r (x, 1)$ and $S_s (x, 1)$ of degree $r$ and $s$ respectively given by the formulas
\begin{equation}\label{def-T-S}
T_r (x, 1) =\sum_{i=1}^r t_i^2(x-\a_i)^2,   \quad \textit{and } \quad S_s (x, 1) =  \sum_{j=1}^s 2u_j^2(x-\b_j)(x-\bar \b_j ), 
\end{equation} 
%
%\begin{equation} \q_f (x, y) =\sum_{i=1}^r t_i^2(x-\a_iy)^2+ \sum_{j=1}^s 2u_j^2(x-\b_jy)(x-\bar \b_j y), \end{equation}
%
where $t_i$, $u_j$ are   to be determined; see \cite{beshaj-2}.    Then, 
\begin{equation}\label{T-S}
\begin{split}
T_r (x, 1) &  =  \left( \sum_{i=1}^r  t_i^2  \right)  x^2  - 2 \left(  \sum_{i-1}^r    t_i^2 \a_i  \right)  x  + \left(  \sum_{i-1}^r    t_i^2 \a_i^2 \right)   \\
S_s (x, 1) & = 2 \left( \sum_{j=1}^s  u_j^2  \right)  x^2 - 4\left(\sum_{j=1}^s     u_j^2  \re (\b_j)  \right)  x + 2 \left(  \sum_{j=1}^s    u_j^2 \cdot ||\b_j||^2 \right).
\end{split}
\end{equation}
For a binary form $f (x, 1)$ of signature $(r, s)$ the quadratic form $\q_f$ is defined as
\begin{equation} \label{Julia-quadratic}
\q_f (x, 1)  :=  T_r (x, 1) + S_s (x, 1) 
\end{equation} 
The discriminant of $\q_f$ is a degree 4 homogenous polynomial in $t_1, \dots t_r,$ $u_1, \dots , u_s$.  We would like to pick values for $t_1, \dots t_r, u_1, \dots , u_s$ such that this discriminant is square free and minimal.  Then we can use the reduction theory of quadratics (with square free, minimal discriminant) to determine the reduced form for $\q_f$. 

For quadratics $T$ and $S$ in Eq.~\eqref{def-T-S} we define 
\begin{equation}
\theta_T = \frac { a_0^2 \cdot \D_T} {t_1^2 \cdots t_r^2}, \qquad \theta_S = \frac {a_0^2 \cdot \D_S} {u_1^4 \cdots u_s^4}
\end{equation}
Notice that  $T_r$ and $S_s$ are given recursively as
\begin{equation}  T_r = T_{r-1} + t_r^2 (x- \a_r ))^2, \qquad S_s = S_{s-1} + u_s^4 \left( x^2 - 2a_s x + (a_s^2+b_s^2)     \right)    
 \end{equation}
Fort $f  \in V_{d, \R}$ with signature $(r, s)$ and equation \cref{deg-n-real-factored},   $\q_f$ is a positive definite quadratic form with discriminant $\Df$  given by the formula
\begin{equation}\label{D-Q}
\begin{split}
\Df = & \D (T_r) + \D (S_s) - 8 \sum_{i, j} t_i^2 u_j^2 \left(  (\a_i-a_j)^2 + b_j^2 \right),
\end{split}
\end{equation}
 see \cite{beshaj-2}.    Let $\theta_0$ of a binary form be 
\begin{equation}
\theta_0 (f) := \frac{a_0^2\cdot |\Df |^{d/2}}{\prod_{i=1}^rt_i^2\, \prod_{j=1}^s u_j^4  }
\end{equation}
and consider $\theta_0(t_1, \dots, t_r, u_1, \dots, u_s)$  as a multivariable function in the variables $t_1, \dots , t_r$, $u_1, \dots , u_s$.  We would like to pick such variables such that $\q_f$ is a reduced quadratic, hence $\Df$ is minimal.  This is equivalent with $\theta_0(t_1, \dots, t_r, u_1, \dots, u_s)$ obtaining a minimal value.   

\begin{lem}\label{unique-min}
The function  $\theta_0 : \R^{r+s} \to \R$  obtains a minimum at a unique point $(\bar t_1, \dots, \bar t_r,  \bar u_1, \dots, \bar u_s)$. 
\end{lem}

Julia in his thesis  \cite{julia} proves existence and Stoll and Cremona prove uniqueness in \cite{stoll-cremona}.   
Choosing $(\bar t_1, \dots, \bar t_r,  \bar u_1, \dots, \bar u_s)$  that   make $\theta_0$ minimal  gives a unique positive definite quadratic $\q_f (x, y)$.  We call this unique quadratic $\q_f (x, y)$ for such a choice of  $(\bar t_1, \dots, \bar t_r,  \bar u_1, \dots, \bar u_s)$  the \textbf{Julia's quadratic} of $f(x, y)$,  denote it by $\J_f (x, y)$, and the quantity  $\theta_f:= \theta_0(\bar t_1, \dots, \bar t_r,  \bar u_1, \dots, \bar u_s)$ the \textbf{ Julia invariant}.  
%From the previous remarks, this is well defined.   

%The following lemma shows that $\theta$ is an invariant of binary forms and $\J$ a covariant of order 2.

\begin{lem}\label{theta-invariant}
%Let $f\in V_{d, \R}$. Then 
Consider $\GL_2 (\C)$ acting on $V_{d, \R}$. Then,    $\theta$ is an   invariant  and   $\J$ is a  covariant  of order 2. 
\end{lem}

%*************************************************
\subsubsection{Minimal height of a binary form in its $\SL_2(\O_k)$-orbit.}   
Consider $f$  a degree $d$ binary form and $k$ its minimal field of definition. Let $M \in \SL_2 (\O_k)$ be a matrix such that $f^M$ is reduced, i.e. $\bar \z(f^M) \in \F_k$ where $\F_k$ is the fundamental domain of $\SL_2(\O_k)$ acting on $\H_3$; see \cite{beshaj} for details.   A  bound on the height of the reduced binary form with respect to Julia invariant is given below. 

\begin{lem}\label{height-bound-theta}
Let   $f(x, 1) = a_0 \prod_{i=1}^d (x-\a_i)$ 
be a reduced binary form where $\a_i$  are the roots. Then, the height of this form   can be bounded by Julia's invariant  as 
\begin{equation} 
H(f) \leq c \cdot \theta_f^{d/2}, \quad \text{ where } \quad  c=    \left( \frac{1}{3}\right)^{\frac{d^2}4}   \,  \left( \frac{4}{d-1}\right)^{\frac{d(d-1)}{2}}    \frac{ 1}{a_0^{d}} \end{equation}
\end{lem}

Let $f$ be a binary form and  $\mathbb F$ its minimal field of definition.   If $f$ is reduced over $\mathbb F$, then it has minimal height in its $\G_{\mathbb F}$-orbit. For    a degree $d$ binary form $f$  defined over $\mathbb F$,
% and $\J_f$ its Julia quadratic,  
$\Df$ its discriminant, and $L  = {\mathbb F} (\Df)$. Then, $[L: {\mathbb F}] \leq d$. 
%If $r$ is the class number of $\J_f$ over $L$ and  $M_1, \dots, M_r$  the matrices with entries in   $\SL_2(\O_k)$ that send $\J_f$ respectively to  $\{ J_1, \dots , J_r\}$, then    the form  $f^{M_j}$ for some $j=1, \dots, r$  has   minimal height over $\SL_2(\O_k)$; see \cite{beshaj} for a proof. 

%***********************************************************************
\subsection{Reduction of binary forms via the hyperbolic centroid}\label{hyp-centroid}
Another reduction method was introduced in \cite{e-sh} following the general idea of \cite{julia}.  The main difference is the definition of the zero map in \cref{zero-map-real}.
We briefly summarize it here. 
Let $V^+_{2n,\R}(0,n)$ denote the set of degree $2n$  binary forms   with real coefficients and no real roots (for the case of binary forms with real roots see \cite{e-sh}). Such binary forms are called \textbf{totally complex}. 
Every $f(x,z)\in V^+_{2n,\R}(0,n)$ can be factored  as  
\begin{equation}
f(x,z)=\prod_{j=1}^n \q_{\alpha_j}(x,z),
\end{equation}
where   $\alpha_j=x_j+{\bf i}y_j$   and  $\q_{\alpha_j}(x,z)=(x-\alpha_jz)(x-\overline{\alpha_j}z)$ are quadratics with real coefficients.
The \textbf{hyperbolic centroid},   ${\mathcal C_{\H}}(\alpha_1,\alpha_2,...,\alpha_n)$ of the collection of a set of points in the upper-half plane  
$\{\alpha_j = x_j + {\bf i} y_j \in \H_2~  \mid  \, j=1,2,...,n\}$ is the unique point $t+{\bf i} u\in \H_2$ that minimizes 
\begin{equation}
\sum_{j=1}^n \frac{(t-x_j)^2+(u-y_j)^2}{uy_j},
\end{equation}
see \cite[Definition~5]{e-sh}.   Let ${\mathfrak s_{j}}$ denote the $j$-th symmetric polynomial in $y_1, \ldots , y_n$.
It follows from (\cite[Prop.~10]{e-sh}) that the centroid  ${\mathcal C_{\H}}=t+{\bf i}u\in \H_2$ of $\alpha_1,\alpha_2,...,\alpha_n$ satisfies
\begin{equation}\label{centroid}
\begin{split}
 \displaystyle{t} & \displaystyle{=\sum_{i=1}^n\left(\frac{ y_1 y_2 \cdots y_{i-1} y_{i+1} \cdots y_n}{\mathfrak s_{n-1}(y_1,y_2,...,y_n)}\right)x_i}  \\
\displaystyle{|{\mathcal C_{\H}}|^2}  & =  \displaystyle{ \sum_{i=1}^n \left(\frac { y_1 y_2 \cdots y_{i-1} y_{i+1} \cdots y_n }   {\mathfrak s_{n-1}(y_1,y_2,...,y_n) }\right)|\alpha_i|^2}  \\
\displaystyle{\q_{{\mathcal C_{\H}}}(x,z)} & = \displaystyle{\sum_{i=1}^n \left(\frac { y_1 y_2 \cdots y_{i-1} y_{i+1} \cdots y_n } {\mathfrak s_{n-1}(y_1,y_2,...,y_n)}\right)\q_{\alpha_i}(x,z)}. 
\end{split}
\end{equation}
The \textbf{centroid zero map} $\z_{\mathcal C}: V^+_{2n,\R}(0,n)\rightarrow \H_2$ is defined via 
\begin{equation}
\z_{\mathcal C}(f):=\mathcal C_{\H}={\mathcal C}_{\H}(\alpha_1,\alpha_2,...,\alpha_n).
\end{equation}
The form
\begin{equation}
\J^{\mathcal C}_f:=  (x-{\mathcal C_{\H}}z)(x-\overline{\mathcal C_{\H}}z)
= \displaystyle{\sum_{j=1}^n \left(\frac { y_1 y_2 \cdots y_{j-1} y_{j+1} \cdots y_n } {\mathfrak s_{n-1}(y_1,y_2,...,y_n)}\right)\q_{\alpha_j}(x,z)}
\end{equation}
is called the \textbf{centroid quadratic} of $f$.
 %\begin{rem}
The reduction theory based on the centroid proceeds as before. Let $f(x,z)$ be a real binary form with no real roots. If $\z_{\mathcal C}(f)\in \F$ then $f$ is reduced. Otherwise, let $M\in \SL_2(\R)$ such that $M^{-1}\z_{\mathcal C}(f) \in \F$. The form $f$ reduces to $f^M(x,z)$.
%\end{rem}
%
\begin{exa}[Totally complex sextics]
Let $f(x,z) \in \Z[x, z] $ be a totally complex sextic factored over $\R$ as 
\[ f(x,z)= (x^2 + a_1 x z + b_1 z^2) (x^2 + a_2 x z + b_2 z^2) (x^2 + a_3 x z + b_3 z^2). \]
Let $d_j = \sqrt{4b_j-a_j^2},~{\bf d}=(d_1, d_2, d_3),~{\bf a} = (a_1, a_2, a_3),~{\bf b}=(b_1,b_2,b_3)$. 
The centroid zero map $\z_{\mathcal C}(f)=t+{\bf i}u\in H_2$ of $f$ is determined by 
\begin{equation}
\begin{split}
t & =- \frac{1}{2} \left( \frac{d_2 d_3}{\mathfrak s_2(d_1,d_2,d_3)}a_1 + \frac {d_1d_3} {\mathfrak s_2(d_1,d_2,d_3)} a_2 +
\frac{d_1 d_2}{\mathfrak s_2(d_1,d_2,d_3)}a_3\right), \\
 |\z_{\mathcal C}(f)|^2 &= \frac{d_2d_3}{\mathfrak s_2(d_1,d_2,d_3)}b_1+\frac{d_1d_3}{\mathfrak s_2(d_1,d-2,d_3)}b_2+\frac{d_1d_2}{\mathfrak s_2(d_1,d_2,d_3)}b_3.
 \end{split}
 \end{equation}
The centroid quadratic of $f$ is given by 
\begin{equation}
\begin{split}
\displaystyle{\J^{\mathcal C}_f }=  &\frac 1 { {\mathfrak s}_2 (d_1,d_2,d_3)} \left(d_2d_3(x^2+a_1xz+b_1z^2)+d_1d_3(x^2+a_2xz+b_2z^2)\right.  \\
& \left. +d_1d_2(x^2+a_3xz+b_3z^2)
\right).
\end{split}
\end{equation}
The reduction is defined over $\Q(d_1, d_2, d_3)$; see \cite[Prop.~11]{e-sh}
\end{exa}
This is generalized   to any degree in \cite[Prop.~12]{e-sh}. 
%
%\begin{prop}
Let $f(x,z)$ be a totally complex form factored over $\R$ as below
\[ f(x,z)= \prod_{i=1}^n (x^2 + a_i x z  + b_i z^2) \]
and  $d_i := \sqrt{4b_i-a_i^2}$, for $i=1, \dots ,  n$. 
Let $\mathfrak s_{n-1} := \sum_{i=1 }^r d_1 \cdots  d_{i-1}\hat{d_i}d_{i+1}\cdots d_r$ where $\hat{d_i}$ denote a missing $d_i$.
%, and 
%\[  {\bf a}=(a_1,...,a_n),   \; \; {\bf b}=(b_1,...,b_n), \; \; {\bf d}=(d_1,...,d_n).\]
%
The centroid quadratic of $f(x,z)$ is given by 
\begin{equation}
\displaystyle{\J^{\mathcal C}_f=\sum_{i=1}^n \left(\frac { d_1 d_2 \cdots d_{i-1} d_{i+1} \cdots d_n } {\mathfrak s_{n-1}}\right)(x^2+a_ixz+b_i z^2)}.
\end{equation}
The centroid zero map $\z_{\mathcal C}(f)=t+{\bf i}u\in \H_2$ is given by
\begin{equation}
\begin{split}
  t & =  - \frac 1 {2}  \,     \sum_{i=1 }^n \frac{d_1 \cdots  d_{i-1}d_{i+1}\cdots d_n}{\mathfrak s_{n-1}}a_i,  \\
   u^2 &  = \frac  1   { 4 \mathfrak s_{n-1}^2}     \, \prod_{i=1}^n d_i      \,  
   \left( \mathfrak s_{n-1} \, \sum_{i=1}^n d_i  +      \sum_{i}^n \, d_1 \cdots \hat{d_i}\cdots \hat{d_j}\cdots d_n  \left( a_i - a_j  \right)^2    \right)    \\
|\z_{\mathcal C}(f)|^2 & =     \sum_{i=1 }^n \frac{d_1 \cdots  d_{i-1}d_{i+1}\cdots d_n}{\mathfrak s_{n-1}}b_i.  \\
%
%\text{where}~\mathfrak s_{n-1} & = \sum_{i=1 }^r d_1 \cdots  d_{i-1}\hat{d_i}d_{i+1}\cdots d_r %~\text{and} ~$\hat{x}$~ \text{denotes a missing}~ $x$ \\   
\end{split}
\end{equation}
The reduction is defined over $\Q(d_1, d_2, ...,d_n)$.
%\end{prop}
%
As pointed out in \cite{e-sh}, expressing  $\z_{\mathcal C}(f)$ in terms of invariants of $f$ or symmetries of the roots of $f$, would be interesting problems on their own.  It is shown in \cite{e-sh} that this hyperbolic reduction is different   from the reduction method in \cite{beshaj, beshaj-2, julia} even though they both seem to find correctly the binary form  with minimal height. 
%In either case 

%***********************************************************
\section{Concluding remarks}
The methods described in this paper are new, as far as we are aware,  and give a new approach for bookkeeping of points in the moduli space of superelliptic curves, provided that we have explicit descriptions of  invariants of binary forms.  For example, for genus $g=2$, since we know explicitly the Igusa arithmetic invariants $J_2, J_4, J_6, J_{10}$ we can explicitly list all points in $\wP_2^3 (k)$ of a given weighted moduli height, including their twists.   This makes it possible to study the arithmetic of the moduli space $\M_2$ and its rational points; see \cite{rat-pts}.       A similar approach can be used for any moduli space of curves when the corresponding invariants are explicitly known.   To our knowledge, the only time when these two types of reduction have been combined is in seminal work of Birch and Swinnerton-Dyer in \cite{birch-s-d-1, birch-s-d-2} for their computation with elliptic curves in an attempt to verify their famous conjecture.   While in this paper we concerned ourselves with a generic superelliptic curve, cases when the moduli point is a singular point or equivalently the curve has a large automorphism group (see \cite{allen, ruben,    kyoto, sh-1})   are even more interesting.

We must point out that the reduction of coefficients is build on the analogy with quadratics and seems to work fine in all the cases. However, they is no known proof as far as we are aware that this reduction will guarantee the binary form with smallest hight in the sense of \cite{heights}.  This seems as a problem worth investigating. 

The reader interested in further details from both geometric and arithmetic aspects of these problems can check \cite{sh-bin-oct, r-sh, milagros, obus-sh, sh-2, curri}.

%\cite{ajm-jones, nebe, ajm-tseng, ajm-gromadzki, obus, ajm-razon,   frey-sh, b-e-sh}
%*************************************************************
%\newpage

% amsplain, unsrt, 

\bibliographystyle{abbrv}
\bibliography{ref-2}
\end{document}